\documentclass[preprint,12pt]{elsarticle}

 \usepackage{graphics}
\usepackage{amsmath}
\usepackage{amssymb}
\usepackage{amsthm}
\usepackage{lineno}

\theoremstyle{plain}

\newtheorem{theorem}{Theorem}[section]
\newtheorem{lemma}[theorem]{Lemma}
\newtheorem{proposition}{Proposition}

\theoremstyle{definition}
\newtheorem{definition}{Definition}
\theoremstyle{example}

\theoremstyle{remark}
\newcommand{\pff}{{\it Proof.}\ }

\begin{document}

\begin{frontmatter}



\title{Shapes of RNA  pseudoknot structures}


\author{Christian M. Reidys$^{*,\dagger}$ and Rita R. Wang $^{*}$}
\address{$*$ Center for Combinatorics, LPMC-TJKLC\\
         $\dagger$ College of Life Sciences    \\
         Nankai University  \\
         Tianjin 300071\\
         P.R.~China\\
         Phone: *86-22-2350-6800\\
         Fax:   *86-22-2350-9272\\
         duck@santafe.edu}

\begin{abstract}
In this paper we study abstract shapes of $k$-noncrossing,
$\sigma$-canonical RNA pseudoknot structures. We consider ${\sf
lv}_k^{\sf 1}$- and ${\sf lv}_k^{\sf 5}$-shapes, which represent a
generalization of the abstract $\pi'$- and $\pi$-shapes of RNA
secondary structures introduced by \citet{Giegerich:04ashape}. Using
a novel approach we compute the generating functions of ${\sf
lv}_k^{\sf 1}$- and ${\sf lv}_k^{\sf 5}$-shapes as well as the
generating functions of all ${\sf lv}_k^{\sf 1}$- and ${\sf
lv}_k^{\sf 5}$-shapes induced by all $k$-noncrossing,
$\sigma$-canonical RNA structures for fixed $n$. By means of
singularity analysis of the generating functions, we derive explicit
asymptotic expressions.
\end{abstract}


\end{frontmatter}

\section{Introduction}\label{S:Intro}

\begin{figure}[ht]
\centerline{\includegraphics[width=0.9\textwidth]{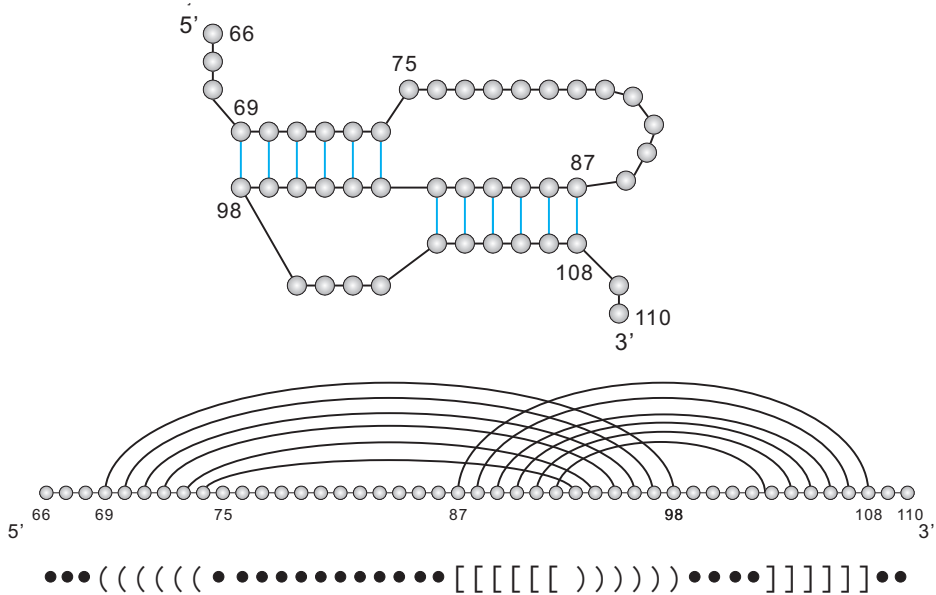}}
\caption{\small The pseudoknot structure of the PrP-encoding mRNA.}
\label{F:pse}
\end{figure}

Pseudoknots have long been known as important structural elements
\cite{Westhof:92a}, see Fig.~\ref{F:pse}. They represent cross-serial
interactions between RNA nucleotides and are an important functionally
in tRNAs, RNaseP \cite{Loria:96a}, telomerase RNA \cite{Butcher},
and ribosomal RNAs \cite{Konings:95a}. Pseudoknots in plant virus RNAs mimic
tRNA structures, and {\it in vitro} selection experiments have produced
pseudoknotted RNA families that bind to the HIV-1 reverse transcriptase
\cite{Tuerk:92}. Import general mechanism, such as ribosomal frame shifting,
are dependent upon pseudoknots \cite{Chamorro:91a}.

Despite their biological importance, pseudoknots are typically
excluded from large-scale computational studies. Although the
problem has attracted considerable attention in the last decade, and
several software tools \cite{Fenix:08,Rivas:99} have become
available, the required resources have remained prohibitive for
applications beyond individual molecules.

\begin{figure}[ht]
\centerline{\includegraphics[width=0.9\textwidth]{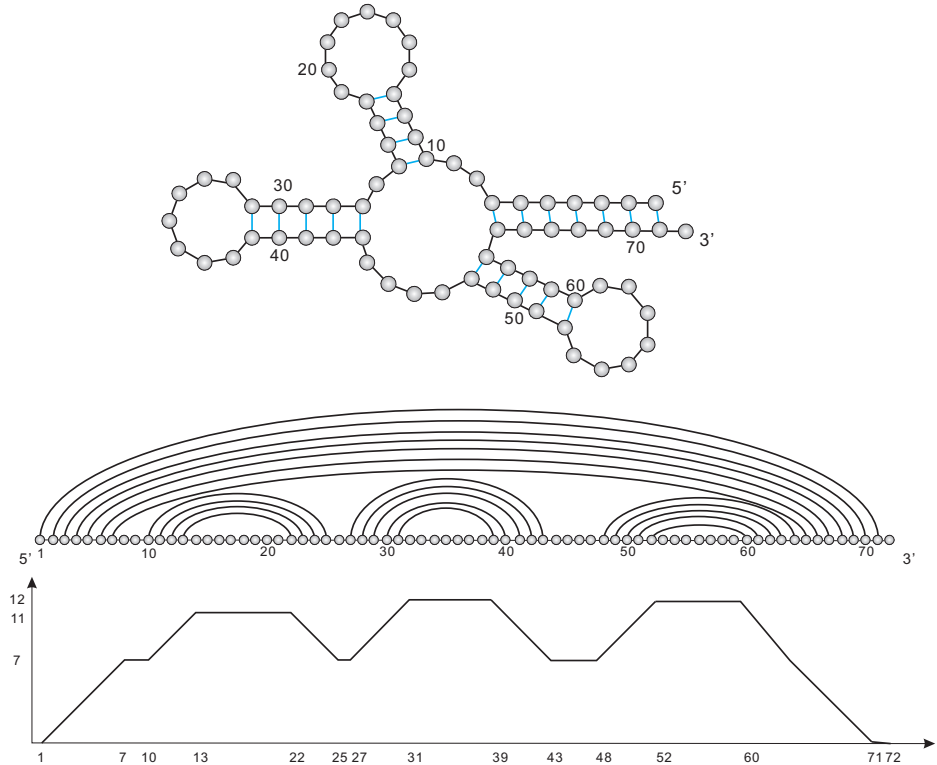}}
\caption{\small { The Sprinzl tRNA RD7550 secondary structure
represented as a planar graph (top), $2$-noncrossing diagram
(middle) and Motzkin-path (bottom), where up/down/horizontal-steps
correspond to start/end/unpaired vertices, respectively.}}
\label{F:Pnas2}
\end{figure}
\begin{figure}[ht]
\centerline{\includegraphics[width=0.85\textwidth]{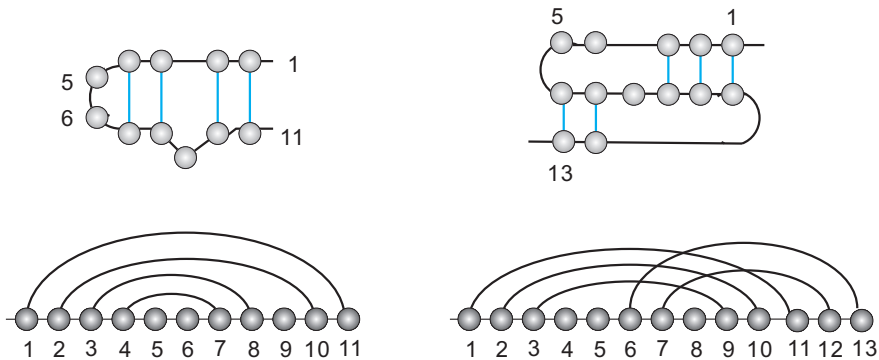}}
\caption{\small A $2$-noncrossing, $2$-canonical RNA structure
(left) and a $3$-noncrossing, $2$-canonical RNA structure (right).}
\label{F:canonical}
\end{figure}

An RNA molecule is a sequence of the four nucleotides {\bf A}, {\bf
G}, {\bf U} and {\bf C} together with the Watson-Crick ({\bf A-U},
{\bf G-C}) and {\bf U-G} base pairing rules. The sequence of bases
is called the primary structure of the RNA molecule. Two bases in
the primary structure which are not adjacent may form hydrogen bonds
following the Watson-Crick base pairing rules. Three decades ago
Waterman {\it et al.} \cite{Kleitman:70,Nussinov,Waterman:79}
analyzed RNA secondary structures. Secondary structures are coarse
grained RNA contact structures. They can be represented as diagrams,
planar graphs as well as Motzkin-paths, see Fig.~\ref{F:Pnas2}.
Diagrams are labeled graphs over the vertex set $[n]=\{1, \dots,
n\}$ with vertex degrees $\le 1$, represented by drawing its
vertices on a horizontal line and its arcs $(i,j)$ ($i<j$), in the
upper half-plane, see Fig.~\ref{F:Pnas2} and Fig.~\ref{F:canonical}.
Here, vertices and arcs correspond to the nucleotides {\bf A}, {\bf
G}, {\bf U} and {\bf C} and Watson-Crick ({\bf A-U}, {\bf G-C}) and
({\bf U-G}) base pairs, respectively. In a diagram two arcs
$(i_1,j_1)$ and $(i_2,j_2)$ are called crossing if $i_1<i_2<j_1<j_2$
holds. Accordingly, a $k$-crossing is a sequence of arcs
$(i_1,j_1),\dots,(i_k,j_k)$ such that
$i_1<i_2<\dots<i_k<j_1<j_2<\dots <j_k$, see Fig.~\ref{F:canonical}.
We call diagrams containing at most $(k-1)$-crossings,
$k$-noncrossing diagrams ($k$-noncrossing partial matchings).

An important observation in this context is that RNA secondary structures have
no crossings in their diagram representation, see Fig.~\ref{F:canonical}
(l.h.s.) and Fig.~\ref{F:Pnas2}, and are therefore $2$-noncrossing diagrams.
 The length of an arc $(i,j)$ is given by $j-i$, characterizing the minimal
length of a hairpin loop.
A stack of length $\sigma$ is a sequence of ``parallel'' arcs of the form
\begin{equation}\label{E:stack}
((i,j),(i+1,j-1),\ldots, (i+(\sigma-1),j-(\sigma-1))).
\end{equation}
In the context of minimum-free energy pseudoknot structures
\cite{Fenix:08} a minimum stack length $\sigma$ or either two or
three is stipulated. We remark that RNA secondary structures are
$2$-noncrossing, $2$-canonical diagrams, whose numbers are
asymptotically given by \cite{Schuster:DAM}
\begin{equation}
S_{2,2}(n) \sim c\, n^{-3/2}\, 1.96798^n, \quad c>0.
\end{equation}
We call an arc of length one a $1$-arc. A $k$-noncrossing, $\sigma$-canonical
RNA structure is a $k$-noncrossing diagram without $1$-arcs, having a minimum
stack-size of $\sigma$.

The efficient minimum free energy (mfe) folding of secondary structures
is a consequence of the following relation of the numbers of RNA secondary
structures over $n$ nucleotides, $S_2(n)$, \cite{Waterman:79}
\begin{equation}\label{E:sn}
S_2(n)=S_2(n-1)+\sum_{j=0}^{n-2}S_2(n-2-j)S_2(j),
\end{equation}
where $S_2(n)=1$ for $0\le n\le 2$. Accordingly, RNA secondary
structures satisfy a constructive recursion. As mentioned above,
this relation is the key for deriving the fundamental DP-recursions
used for the polynomial time folding of secondary structures
\cite{Vienna:Server,Nussinov} and has therefore profound algorithmic
implications. In addition, eq.~(\ref{E:sn}) is of central importance
for the analysis of abstract shapes \cite{Nebel}. In addition, for a
given RNA sequence, we have not only one but an ensemble of
structures, quantified via the partition function generated by the
(Boltzman weighted) probability space of all structures
\cite{McCaskill}. In view of the fact that the number of the mfe and
suboptimal foldings of an RNA sequence is large, Giegerich {\it et
al.} \cite{Giegerich:04ashape} introduced the notion of abstract
shapes of secondary structures. Two particularly important shape
levels are the important level-1 ($\pi'$-) and level-5 ($\pi$-)
shapes were studied in \cite{Giegerich:04ashape}. In
\cite{Giegerich:06shape}, the authors compute the probability of a
shape by means of the partition function, where the probability of a
shape is the induced probability of all the structures inducing it.

\begin{figure}[ht]
\centerline{\includegraphics[width=0.9\textwidth]{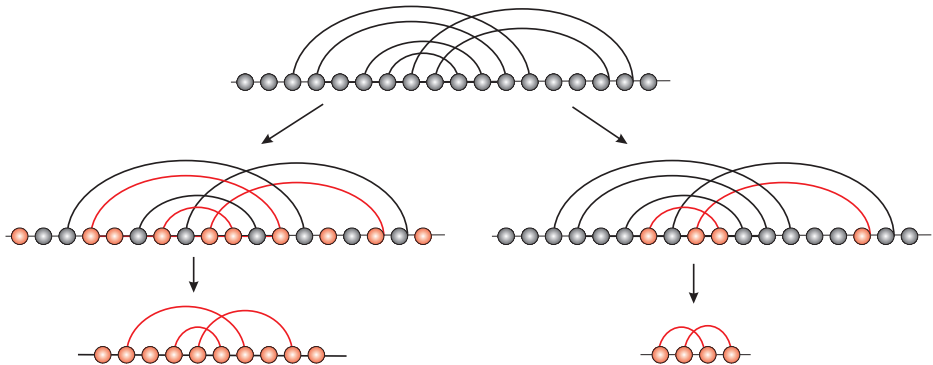}}
\caption{\small ${\sf lv}_k^{\sf 1}$- and ${\sf lv}_k^{\sf
5}$-shapes: a $3$-noncrossing, $2$-canonical RNA structure (top),
its ${\sf lv}_3^{\sf 1}$-shape (bottom left) and its ${\sf
lv}_3^{\sf 5}$-shape (bottom right).} \label{F:shapeintro}
\end{figure}

The problem with pseudoknotted structures is, that they do not
satisfy a recursion of the type of eq.~(\ref{E:sn}), rendering the
{\it ab initio} folding into mfe configurations
\cite{Fenix:08,Lyngso} as well as the derivation of any other
properties a nontrivial task. Here, we generalize the $\pi'$- and
$\pi$-shapes of \cite{Giegerich:04ashape}, by introducing ${\sf
lv}_k^{\sf 1}$- and ${\sf lv}_k^{\sf 5}$-shapes, see
Fig.~\ref{F:shapeintro}.
Our results are not new in case of $k=2$, since we have ${\sf
lv}_2^{\sf 1}=\pi'$ and ${\sf lv}_2^{\sf 5}=\pi$. In two beautiful
papers \cite{Lorenz,Nebel} $\pi'$- and $\pi$-shapes have been
analyzed. The results of \cite{Lorenz,Nebel} explicitly make use of
the constructive recurrence relation given in eq.~(\ref{E:sn}).
Their approach can consequently not be generalized to RNA pseudoknot
structures, as the latter are genuinely nonrecursive. Our framework
therefore identifies the combinatorial ``heart'' of the results of
\cite{Lorenz,Nebel} and provides a new approach avoiding any notion
of grammar or recursiveness. The key idea behind the construction of
${\sf lv}_k^{\sf 1}$- and ${\sf lv}_k^{\sf 5}$-shapes is a
projection onto so called $k$-noncrossing core-structures
\cite{Reidys:07lego}.

The paper is organized as follows: after introducing all necessary
background we give a detailed computation of the generating
functions and study their singularities. We derive simple asymptotic
expressions for the numbers of ${\sf lv}_k^{\sf 1}$- and ${\sf
lv}_k^{\sf 5}$-shapes as well as the numbers of theses shapes,
induced by $k$-noncrossing, $\sigma$-canonical RNA structures of
fixed length $n$. Finally we put our results into context.

\section{Some basic facts}\label{S:Preli}

Let $f_k(n,\ell)$ denote the number of $k$-noncrossing diagrams
on $n$ vertices having exactly $\ell$ isolated vertices. A diagram
without isolated points is called a matching. The exponential generating
function of $k$-noncrossing matchings satisfies the following identity
\cite{Chen,Grabiner:93a,Reidys:07pseu}
\begin{eqnarray}\label{E:ww0}
\label{E:ww1} \sum_{n\ge 0} f_{k}(2n,0)\cdot\frac{z^{2n}}{(2n)!} & =
& \det[I_{i-j}(2z)-I_{i+j}(2z)]|_{i,j=1}^{k-1}
\end{eqnarray}
where $I_{r}(2z)=\sum_{j \ge 0}\frac{z^{2j+r}}{{j!(j+r)!}}$ is the
hyperbolic Bessel function of the first kind of order $r$.
Eq.~(\ref{E:ww0}) allows to conclude that the ordinary generating
function
\begin{equation*}
{\bf F}_k(z)=\sum_{n\ge 0}f_k(2n,0) z^{n}
\end{equation*}
is $D$-finite \cite{Stanley:80}, i.e.\ there exists some $e\in \mathbb{N}$
such that
\begin{equation}
q_{0,k}(z)\frac{d^e}{d z^e}{\bf F}_k(z)+q_{1,k}(z)\frac{d^{e-1}}{d
z^{e-1}}{\bf F}_k(z)+\cdots+q_{e,k}(z){\bf F}_k(z)=0,
\end{equation}
where $q_{j,k}(z)$ are polynomials. Since $I_{r}(2z)$ is $D$-finite
by its definition and $D$-finite power series are algebraic closed
\cite{Stanley:80}. The key point is that any singularity of ${\bf F}_k(z)$
is contained in the set of roots of $q_{0,k}(z)$ \cite{Stanley:80}, which
we denote by $R_k$.
For $2\leq k\leq 9$, we give the polynomials $q_{0,k}(z)$ and their roots in
Table \ref{Table:polyroot}.

\begin{table}
\begin{center}
\begin{tabular}{cll}
\hline
$k$ & $q_{0,k}(z)$ & $R_k$  \\
\hline
$2$ & $(4z-1)z$ & $\{\frac{1}{4}\}$\\
$3$ & $(16z-1)z^2$ & $\{\frac{1}{16}\}$\\
$4$ & $(144z^2-40z+1)z^3$ & $\{\frac{1}{4},\frac{1}{36}\}$\\
$5$ & $(1024z^2-80z+1)z^4$ & $\{\frac{1}{16},\frac{1}{64}\}$\\
$6$ & $(14400z^3-4144z^2+140z-1)z^5$ & $\{\frac{1}{4},\frac{1}{36},\frac{1}{100}\}$\\
$7$ & $(147456z^3-12544z^2+224z-1)z^6$ & $\{\frac{1}{16},\frac{1}{64},\frac{1}{144}\}$\\
$8$ & $(2822400z^4-826624z^3+31584z^2-336z+1)z^7$ & $\{\frac{1}{4},\frac{1}{36},\frac{1}{100},\frac{1}{196}\}$\\
$9$ & $(37748736z^4-3358720z^3+69888z^2-480z+1)z^8$ & $\{\frac{1}{16},\frac{1}{64},\frac{1}{144},\frac{1}{256}\}$\\
\hline
\end{tabular}
\centerline{}
\smallskip
\caption{\small We present the polynomials $q_{0,k}(z)$  and their
nonzero roots obtained by the MAPLE package {\tt GFUN}. }
\label{Table:polyroot}
\end{center}
\end{table}
In \cite{Reidys:08k} we showed that for arbitrary $k$
\begin{equation}\label{E:theorem}
f_{k}(2n,0) \, \sim  \, \widetilde{c}_k  \, n^{-((k-1)^2+(k-1)/2)}\,
(2(k-1))^{2n},\qquad \widetilde{c}_k>0 \ .
\end{equation}
in accordance with the fact that ${\bf F}_k(z)$ has the unique dominant
singularity $\rho_k^2$, where $\rho_k=1/(2k-2)$.

Let $\mathcal{T}_{k,\sigma}(n)$ denote the set of $k$-noncrossing,
$\sigma$-canonical RNA structures of length $n$ and let
$\mathrm{T}_{k,\sigma}(n)$ denote their number.
$\mathcal{T}_{k,\sigma}(n)$ can be identified with the set of
$k$-noncrossing RNA structures with each stack size $\geq \sigma$.
Furthermore, let $\mathcal{T}_{k,\sigma}(n,h)$ denote the set of
$k$-noncrossing, $\sigma$-canonical RNA structures of length $n$
with $h$ arcs, and set
$\mathrm{T}_{k,\sigma}(n,h)=|\mathcal{T}_{k,\sigma}(n,h)|$. The
bivariate generating function of $\mathrm{T}_{k,1}(n,h)$ $(k\geq 2)$
has been computed in \cite{Reidys:07asym}
\begin{equation}
\sum_{n\geq 0}\sum_{h=0}^{\lfloor\frac{n}{2}\rfloor}
\mathrm{T}_{k,1}(n,h)v^{h}y^n=
\frac{1}{vy^2-y+1}{\bf F}_k\left(\frac{vy^2}{\left(vy^2-y+1\right)^2}\right)
\end{equation}
and the generating function for $k$-noncrossing, $\sigma$-canonical
RNA structures is given by \cite{Reidys:07lego}
\begin{equation}\label{E:f1}
\sum_{n \ge 0}{\rm T}_{k,\sigma}^{}(n)y^n =
\frac{1}{u_0y^2-y+1}{\bf F}_k\left(
\left(\frac{\sqrt{u_0}y}{u_0y^2-y+1}\right)^2\right) \ ,
\end{equation}
where $u_0=\frac{(y^2)^{\sigma-1}}{(y^2)^{\sigma}-y^2+1}$.

According to Pringsheim's Theorem \cite{Flajolet:07a,Tichmarsh:39Pringsheim},
each power series $f(z)=\sum_{n\geq 0}a_n\, z^n$
with nonnegative coefficients and a radius of convergence $R>0$ has a
positive real dominant singularity at $z=R$. This singularity plays
a key role for the asymptotics of the coefficients.
The class of theorems that deal with such deductions are called
transfer-theorems \cite{Flajolet:07a}. One key ingredient in this
framework is a specific domain in which the functions in question
are analytic, which is ``slightly'' bigger than their respective
radius of convergence. It is tailored for extracting the
coefficients via Cauchy's integral formula. Details on the method
can be found in \cite{Flajolet:07a,Stanley:80}. In case of
$D$-finite functions we have analytic continuation in any simply
connected domain containing zero \cite{Wasow:87} and all
prerequisites of singularity analysis are met. We use the notation
\begin{equation}\label{E:genau}
\left\{f(z)=O\left(g(z)\right) \ \text{\rm as $z\rightarrow
\rho$}\right\}\quad \Longleftrightarrow \quad
\left\{\frac{f(z)}{g(z)} \ \text{\rm is bounded as $z\rightarrow
\rho$}\right\}.
\end{equation}
Let $[z^n]f(z)$ denote the $n$-th coefficient of the power series
$f(z)$ at $z=0$.
\begin{theorem}\label{T:transfer1}\cite{Flajolet:07a}
Let $f(z),g(z)$ be $D$-finite functions with unique dominant
singularity $\rho$ and suppose
\begin{equation}
f(z) = O( g(z)) \ \mbox{ as }\, z\rightarrow \rho \ .
\end{equation}
Then we have
\begin{equation}
[z^n]f(z)= \,C \,\left(1-O(\frac{1}{n})\right)\,  [z^n]g(z)
\end{equation}
where $C$ is a constant.
\end{theorem}
Theorem \ref{T:transfer1} implies the following result, tailored for
our functional equations. It is a particular instance of the
supercritical paradigm, where we have the following situation: we
are given a $D$-finite function, $f(z)$ and an algebraic function
$g(u)$ satisfying $g(0)=0$. Furthermore we suppose that $f(g(u))$
has the unique real valued dominant singularity $\gamma$ and $g$ is
regular in a disc with radius slightly larger than $\gamma$. The
supercritical paradigm then stipulates that the subexponential
factors of $f(g(u))$ at $u=0$ coincide with those of $f(z)$.

\begin{proposition}\label{P:algeasym}
Suppose $\vartheta_\sigma(z)$ is an algebraic function, analytic for
$\vert z\vert <\delta$ and satisfies $\vartheta_{\sigma}(0)=0$. Suppose further
$\gamma_{k,\sigma}<\delta$ is the real unique dominant singularity of ${\bf
F}_k(\vartheta_{\sigma}(z))$ and satisfies
$\vartheta_{\sigma}(\gamma_{k,\sigma})=\rho_k^2$. Then
\begin{equation}
[z^n]\,{\bf
F}_k(\vartheta_{\sigma}(z)) \sim c_k  \, n^{-((k-1)^2+(k-1)/2)}\,
\left(\gamma_{k,\sigma}^{-1}\right)^n .
\end{equation}
\end{proposition}

Let ${\mathcal G}_k(n,m)$ denote the set of the
$k$-noncrossing matchings of length $2n$ with $m$ $1$-arcs.
In our first lemma, we will compute the bivariate generating
function of $g_k(n,m)$, i.e.~the number of $k$-noncrossing matchings
of length $2n$ with exactly $m$ $1$-arcs.
\begin{lemma}\label{L:fhairpin}
Suppose $k,n,m\in \mathbb{N}$, $k\geq 2$, $0\leq m\leq n$.
Then $g_k(n,m)$ satisfies the recursion
\begin{equation}
(m+1)g_k(n+1,m+1)=(m+1)g_k(n,m+1)+(2n+1-m)g_k(n,m).
\end{equation}
Furthermore, the generating function
${\bf G}_k(x,y)=\sum_{n\geq0}\sum_{m=0}^{n}g_k(n,m)x^ny^m$ is given by
\begin{equation}
{\bf G}_k(x,y)=\frac{1}{x+1-yx}
{\bf F}_k\left(\frac{x}{(x+1-yx)^2}\right).
\end{equation}
\end{lemma}
\pff Choose a $k$-noncrossing matching $\delta\in {\mathcal
G}_k(n+1,m+1)$ and label one $1$-arc. We have $(m+1)g_k(n+1,m+1)$
different such labeled $k$-noncrossing matchings. On the other hand,
in order to obtain such a labeled matching, we can also insert one
labeled $1$-arc in a $k$-noncrossing matching $\delta'\in {\mathcal
G}_k(n,m+1)$. In this case, we can only put it inside one original
$1$-arc in $\delta'$ in order to preserve the number of the
$1$-arcs. We may also insert a labeled $1$-arc in a $k$-noncrossing
matching $\delta''\in {\mathcal G}_k(n,m)$. In this case, we can
only insert the $1$-arc between two vertices not forming a $1$-arc.
Therefore, we arrive at $(m+1)g_k(n,m+1)+(2n+1-m)g_k(n,m)$
different such labeled matchings and
\begin{equation}\label{E:regk}
(m+1)g_k(n+1,m+1)=(m+1)g_k(n,m+1)+(2n+1-m)g_k(n,m).
\end{equation}
This recursion implies the following partial differential equation
for the generating function
\begin{equation}\label{E:diffgf}
x^{-1}\frac{\partial {\bf G}_k(x,y)}{\partial y}=\frac{\partial
{\bf G}_k(x,y)}{\partial y}+2x\frac{\partial {\bf G}_k(x,y)}{\partial
x}+{\bf G}_k(x,y)-y\frac{\partial {\bf G}_k(x,y)}{\partial y},
\end{equation}
whose general solution is given by
\begin{equation}
{\bf G}_k(x,y)=
\frac{F\left(\frac{yx-1-x}{\sqrt{x}}\right)}{\sqrt{x}},
\end{equation}
where $F(z)$ is an arbitrary function. By definition, we have
$\sum_{m=0}^ng_k(n,m)=f_k(2n,0)$ and
\begin{equation}\label{E:initial}
{\bf G}_k(x,1)=\sum_{n\geq 0}f_k(2n,0)x^n.
\end{equation}
Using eq.~(\ref{E:diffgf}) and eq.~(\ref{E:initial}) we derive
\begin{equation}
{\bf G}_k(x,y)=\frac{1}{x+1-yx} \sum_{n\geq
0}f_k(2n,0)\left(\frac{x}{(x+1-yx)^2}\right)^n,
\end{equation}
whence the lemma. \qed

\section{Combinatorics of ${\sf lv}_k^{\sf 5}$-shapes}\label{S:Comb}

\begin{figure}[ht]
\centerline{\includegraphics[width=1\textwidth]{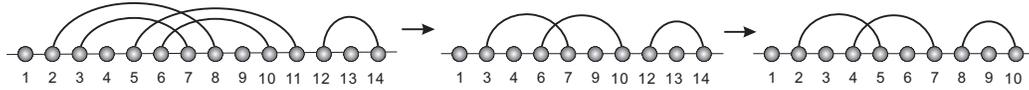}}
\caption{\small A $3$-noncrossing core structure is obtained from a
$3$-noncrossing, $1$-canonical RNA structure in two steps.}
\label{F:core}
\end{figure}

\begin{figure}[ht]
\centerline{\includegraphics[width=0.7\textwidth]{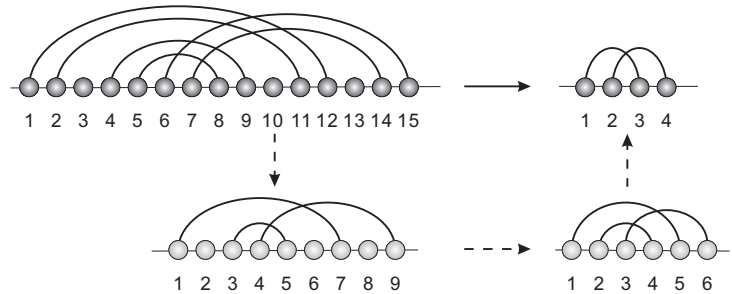}}
\caption{\small Two methods for generating the ${\sf lv}_3^{\sf
5}$-shape. A $3$-noncrossing, $2$-canonical RNA structure (top-left)
is mapped in two ways into its ${\sf lv}_3^{\sf 5}$-shape
(top-right). } \label{F:Ikmapchain}
\end{figure}

We now show how to derive the ${\sf lv}_k^{\sf 5}$-shape of a given
$k$-noncrossing, $\sigma$-canonical RNA structures. This construction
is based on the notion of $k$-noncrossing cores \cite{Reidys:07lego}. A
$k$-noncrossing core is a $k$-noncrossing RNA structure in which each
stack has size exactly one. The cores of a $k$-noncrossing,
$\sigma$-canonical RNA structure, $\delta$, denoted by $c(\delta)$
is obtained in two steps: first we map arcs and isolated vertices
as follows:
\begin{equation}\label{E:coremap}
\forall \ell \geq \sigma-1; \quad ((i-\ell, j+\ell),\ldots, (i,j))
\mapsto (i,j)\text{ and } j\mapsto j\text{ if $j$ is isolated}
\end{equation}
and second we relabel the vertices of the resulting diagram from left
to right in increasing order, see Fig.\ref{F:core}.
We are now in position to define ${\sf lv}_k^{\sf 5}$-shapes.
\begin{definition}{\bf (${\sf lv}_k^{\sf 5}$-shape)}
Given a $k$-noncrossing, $\sigma$-canonical RNA structure $\delta$, its
${\sf lv}_k^{\sf 5}$-shape, ${\sf lv}_k^{\sf 5}(\delta)$, is obtained by first
removing all isolated vertices and second apply the core-map $c$.
\end{definition}
Alternatively the ${\sf lv}_k^{\sf 5}$-shape can also be derived as follows: we
first project into the core $c(\delta)$, second, we remove all
isolated vertices and third we apply the core-map $c$ again,
see Fig.\ref{F:Ikmapchain}.
The second step is a projection from $k$-noncrossing cores to
$k$-noncrossing matchings and surjective, since for each $k$-noncrossing
matching $\alpha$, we can obtain a core structure by inserting isolated
vertices between any two arcs contained in some stack. By construction,
${\sf lv}_k^{\sf 5}$ shapes do not preserve stack-lengths, interior loops
and unpaired regions.

Let ${\mathcal I}_k(n,m)$ ($i_k(n,m)$) denote the set (number) of
the ${\sf lv}_k^{\sf 5}$-shapes of length $2n$ with $m$ $1$-arcs and
\begin{equation}
{\bf I}_k(z,u)=\sum_{n\geq0}\sum_{m=0}^{n} i_k(n,m)z^nu^m
\end{equation}
be the bivariate generating function. Furthermore, let $i_k(n)$ denote the
number of the ${\sf lv}_k^{\sf 5}$-shapes of length $2n$ with generating function
\begin{equation}
{\bf I}_k(z)=\sum_{n\geq0}i_k(n)z^n.
\end{equation}
Since any ${\sf lv}_k^{\sf 5}$-shape is in particular the core of some
$k$-noncrossing matching, Lemma \ref{L:fhairpin} allows us to establish a
relation between the bivariate generating function of $i_k(n,m)$ and the
generating function of ${\bf F}_k(z)$.
\begin{theorem}\label{T:gfIk}
Let $k,n,m$ be natural numbers where $k\geq 2$, then the following assertions
hold\\
{\rm (a)} the generating functions ${\bf I}_k(z,u)$ and ${\bf I}_k(z)$ satisfy
\begin{eqnarray}\label{E:gfIk}
{\bf I}_k(z,u) & = & \frac{1+z}{1+2z-zu}
{\bf F}_k\left(\frac{z(1+z)}{(1+2z-zu)^2}\right)\\
\label{E:o1}
{\bf I}_k(z)   & = & {\bf F}_k\left(\frac{z}{1+z}\right).
\end{eqnarray}
{\rm (b)} for $2\leq k\leq 9$, the number of
${\sf lv}_k^{\sf 5}$-shapes of length $2n$ is
asymptotically given by
\begin{equation}\label{E:asy-p}
i_k(n)\sim c_k n^{-((k-1)^2+(k-1)/2)}\left(\mu_{k}^{-1}\right)^{n},
\end{equation}
where $\mu_{k}$ is the unique minimum positive real solution of
$\frac{z}{1+z}=\rho_k^2$ and $c_{k}$ is some positive constant.
\end{theorem}
\pff
We first prove (a). For this purpose we
define a map between $k$-noncrossing matchings with $m$ $1$-arcs and
${\sf lv}_k^{\sf 5}$-shapes
\begin{equation*}
g\colon \ {\mathcal G}_k(n,m)\rightarrow \dot{\bigcup_{0 \leq b \leq
n-m}}\left[{\mathcal I}_k(n-b,m)\times \left\{(a_j)_{1\leq j\leq
n-b}\mid\sum_{j=1}^{n-b}a_j=b,\ a_j\geq 0\right\}\right],
\end{equation*}
where $n\geq 1$. Here, for every $\delta\in
\mathcal{G}_k(n,m)$, we have $g(\delta)=(c(\delta),(a_j)_{1\leq j\leq
n-b})$, where $c(\delta)$ is the core structure of $\delta$ obtained
according to eq.~(\ref{E:coremap}) and where $(a_j)_{1\leq j\leq
n-b}$ keeps track of the deleted arcs. It is straightforward to
check that the map $g$ is well defined, since all the $1$-arcs of
$c(\delta)$ are just the $1$-arcs of $\delta$. By construction, $g$
is a bijection and we have
\begin{equation*}
|\{(a_j)_{1\leq j\leq n-b}\mid\sum_{j=1}^{n-b}a_j=b,\ a_j\geq
0\}|=\binom{n-1}{b}.
\end{equation*}
Then we derive
\begin{equation}
g_k(n,m)=\sum_{b=0}^{n-m}\binom{n-1}{b}i_k(n-b,m),\quad\text{for
} n\geq 1,
\end{equation}
which implies
\begin{equation*}
\sum_{n\geq 0}
\sum_{m=0}^{n}g_k(n,m)x^ny^m=\sum_{n\geq1}\sum_{m=0}^{n}\sum_{b=0}^{
n-m}\binom{n-1}{b}i_k(n-b,m)x^ny^m+1.
\end{equation*}
We next observe
\begin{equation*}
\sum_{n\geq1}\sum_{m=0}^{n}\sum_{b=0}^{
n-m}\binom{n-1}{b}i_k(n-b,m)x^ny^m =\sum_{b\geq
0}\sum_{m\geq 0}\sum_{n\geq
n_0}\binom{n-1}{b}i_k(n-b,m)x^ny^m,
\end{equation*}
where $n_0=\max\{ m+b,1\}$ and setting $s=n-b$,
\begin{equation*}
\sum_{n\geq1}\sum_{m=0}^{n}\sum_{b=0}^{
n-m}\binom{n-1}{b}i_k(n-b,m)x^ny^m=\sum_{b\geq
0}\sum_{m\geq 0}\sum_{s\geq
s_0}\binom{s+b-1}{b}i_k(s,m)x^{s+b}y^m,
\end{equation*}
where $s_0=\max\{m,1\}$. In view of
\begin{equation*}
\sum_{b\geq 0}\binom{s+b-1}{b}x^b=\frac{1}{(1-x)^s}
\end{equation*}
and interchanging the terms of summation, we derive
\begin{equation*}
\sum_{n\geq1}\sum_{m=0}^{n}\sum_{b=0}^{
n-m}\binom{n-1}{b}i_k(n-b,m)x^ny^m=\sum_{s\geq 1}\sum_{m=
0}^{s}i_k(s,m)\left(\frac{x}{1-x}\right)^sy^m
\end{equation*}
and arrive at
\begin{equation*}
\sum_{n\geq
0}\sum_{m=0}^{n}g_k(n,m)x^ny^m=\sum_{n\geq0}\sum_{m=0}^{n}
i_k(n,m)\left(\frac{x}{1-x}\right)^ny^m.
\end{equation*}
According to Lemma \ref{L:fhairpin}, we have
\begin{equation*}
\sum_{n\geq 0}\sum_{m=0}^{n}g_k(n,m)x^ny^m=\frac{1}{x+1-yx}
\sum_{n\geq 0}f_k(2n,0)\left(\frac{x}{(x+1-yx)^2}\right)^n,
\end{equation*}
setting $z=\frac{x}{1-x}$ and $u=y$,
\begin{equation*}
\sum_{n\geq0}\sum_{m=0}^{n}
i_k(n,m)z^nu^m=\frac{1+z}{1+2z-zu}\sum_{n\geq
0}f_k(2n,0)\left(\frac{z(1+z)}{(1+2z-zu)^2}\right)^n.
\end{equation*}
In particular, setting $u=1$, we derive
\begin{equation*}
\sum_{n\geq0}i_k(n)z^n=\sum_{n\geq
0}f_k(2n,0)\left(\frac{z}{1+z}\right)^n,
\end{equation*}
whence (a) follows.\\
Assertion (b) is a direct consequence of the supercritical paradigm, see
Proposition~\ref{P:algeasym}. As mentioned before, the ordinary generating
function ${\bf F}_k(z)=\sum_{n\ge 0}f_k(2n,0) z^{n}$ is $D$-finite
\cite{Stanley:80} and the inner function $\vartheta(z)=\frac{z}{1+z}$ is
algebraic, satisfies $\vartheta(0)=0$ and is analytic for $\vert z\vert<1$.
By direct calculation, using the fact that all singularities of ${\bf F}_k(z)$
are contained within the set of zeros of $q_{0,k}(z)$, see
Tab.~\ref{Table:polyroot}, we can then verify that ${\bf F}_k(\vartheta(z))$
has the unique dominant real singularity $\mu_k<1$ satisfying
$\vartheta(\mu_k)=\rho_k^2$ for $2\le k\le 9$. In view of
$f_{k}(2n,0) \sim  \widetilde{c}_k  n^{-((k-1)^2+(k-1)/2)}\,(2(k-1))^{2n}$,
Proposition~\ref{P:algeasym} guarantees eq.~(\ref{E:asy-p})
\begin{equation*}
i_k(n)\sim c_k n^{-((k-1)^2+(k-1)/2)}\left(\mu_{k}^{-1}\right)^{n}.
\end{equation*}
This proves (b) completing the proof of the theorem.
\qed

We next studying the number of ${\sf lv}_k^{\sf 5}$-shapes induced by
$k$-noncrossing, $\sigma$-canonical RNA structures of fixed length $n$,
${\sf lv}_{k,\sigma}^{\sf 5}(n)$, setting
\begin{equation}\label{E:lv5}
{\bf Lv}_{k,\sigma}^{\sf 5}(x)=
\sum_{n\geq 0}{\sf lv}_{k,\sigma}^{\sf 5}(n) x^n.
\end{equation}

\begin{theorem}\label{T:Ikasy}
Let $k,\sigma\in \mathbb{N}$, where $k\ge 2$. Then the following assertions hold\\
{\rm (a)} the generating function
${\bf Lv}_{k,\sigma}^{\sf 5}(x)$ is given by
\begin{equation}
{\bf Lv}_{k,\sigma}^{\sf 5}(x)=\frac{(1+x^{2\sigma})}{(1-x)(
1+2x^{2\sigma}-x^{2\sigma+1})}{\bf F}_k\left(\frac{x^{2\sigma}
(1+x^{2\sigma})}{\left(1+2x^{2\sigma}-x^{2\sigma+1}\right)^2}
\right).
\end{equation}
{\rm (b)} for $2\leq k\leq 9$ and $1\leq \sigma \leq 10$
\begin{equation}
{\sf lv}_{k,\sigma}^{\sf 5}(n)\sim
c_{k,\sigma}n^{-((k-1)^2+(k-1)/2)}\left(\zeta_{k,\sigma}^{-1}\right)^{n},
\end{equation}
where $c_{k,\sigma}>0$ and $\zeta_{k,\sigma}$ is the unique minimum positive
real solution of
\begin{equation}\label{E:Ikinner}
\frac{x^{2\sigma}
(1+x^{2\sigma})}{\left(1+2x^{2\sigma}-x^{2\sigma+1}\right)^2}=\rho_k^2.
\end{equation}
\end{theorem}
\begin{table}[ht]
\begin{center}
\begin{tabular}{cccccccc}
\hline
$\sigma/k$& $2$ & $3$ & $4$ & $5$ & $6$ & $7$ & $8$ \\
\hline
$1$ & 1.51243& 3.67528 & 5.77291& 7.82581 & 9.85873 & 11.88118 & 13.89746\\
$2$ &1.26585&1.93496&2.41152&2.80275&3.14338&3.44943&3.72983\\
$3$ & 1.17928& 1.55752& 1.80082 &1.98945&2.14693&2.28376&2.40567\\
\hline
\end{tabular}
\centerline{}
\smallskip
\caption{\small The exponential growth rates $\zeta_{k,\sigma}^{-1}$
of ${\sf lv}_k^{\sf 5}$-shapes induced by $k$-noncrossing, $\sigma$-canonical
RNA structures of length $n$.}
\label{Table:Ik}
\end{center}
\end{table}
\pff In order to proof of (a) we observe that we can always inflate a structure
by adding arcs to stacks or isolated vertices without changing its
${\sf lv}_k^{\sf 5}$-shape.
In fact, for any given ${\sf lv}_k^{\sf 5}$-shape, $\beta$, adding the minimal
number of arcs
to each stack such that every stack has $\sigma$ arcs, and inserting one isolated
vertex in any $1$-arc, we derive a $k$-noncrossing, $\sigma$-canonical structure
having arc-length$\geq 2$, of minimal length.
We can therefore derive ${\bf Lv}_{k,\sigma}^{\sf 5}(x)$, see eq.(\ref{E:lv5}),
from the bivariate generating function ${\bf I}_k(z,u)$ as follows
\begin{equation*}
{\bf Lv}_{k,\sigma}^{\sf 5}(x)=\sum_{n\geq
0}\sum_{s=0}^{\lfloor\frac{n}{2\sigma}
\rfloor}\sum_{m=0}^{\min\{s,n-2\sigma s\}}i_k(s,m)x^{n} =\sum_{s\geq
0}\sum_{m=0}^{s}\sum_{n\geq 2\sigma s+m}i_k(s,m)x^{n},
\end{equation*}
whence
\begin{equation*}
{\bf Lv}_{k,\sigma}^{\sf 5}(x)=\frac{1}{1-x}\sum_{s\geq
0}\sum_{m=0}^{s} i_k(s,m)x^{2\sigma s+m}
\end{equation*}
and in view of eq.~(\ref{E:gfIk}), ${\bf I}_k(z,u)=\frac{1+z}{1+2z-zu}
{\bf F}_k\left(\frac{z(1+z)}{(1+2z-zu)^2}\right)$, we derive
\begin{equation*}
{\bf Lv}_{k,\sigma}^{\sf 5}(x)=\frac{(1+x^{2\sigma})}{(1-x)(
1+2x^{2\sigma}-x^{2\sigma+1})}{\bf F}_k\left(\frac{x^{2\sigma}
(1+x^{2\sigma})}{\left(1+2x^{2\sigma}-x^{2\sigma+1}\right)^2}
\right).
\end{equation*}
As for (b), we observe that the factor
\begin{equation*}
\varphi_\sigma(x)=\frac{(1+x^{2\sigma})}{(1-x)(
1+2x^{2\sigma}-x^{2\sigma+1})}
\end{equation*}
does not induce a dominant singularity of ${\bf Lv}_{k,\sigma}^{\sf 5}(x)$.
Therefore all dominant singularities of ${\bf Lv}_{k,\sigma}^{\sf 5}(x)$ stem from
${\bf F}_k\left(\frac{x^{2\sigma}
(1+x^{2\sigma})}{\left(1+2x^{2\sigma}-x^{2\sigma+1}\right)^2}
\right)$. Indeed, assume {\it a contrario} that there were some dominant
singularity of ${\bf Lv}_{k,\sigma}^{\sf 5}(x)$, $\zeta$, that is induced
by $\varphi_\sigma(x)$. This would imply that
$\zeta$ is also a dominant singularity of ${\bf F}_k\left(\frac{x^{2\sigma}
(1+x^{2\sigma})}{\left(1+2x^{2\sigma}-x^{2\sigma+1}\right)^2}
\right)$ which immediately leads to a contradiction.\\
We next verify that for $2\leq k\leq 9$ and $1\leq \sigma \leq 10$, the
minimum positive real solution of eq.~(\ref{E:Ikinner}), $\zeta_{k,\sigma}$,
is the unique dominant singularity of ${\bf Lv}_{k,\sigma}^{\sf 5}(x)$ and
Proposition~\ref{P:algeasym} implies
\begin{equation*}
{\sf lv}_{k,\sigma}^{\sf 5}(n)\sim
c_{k,\sigma}n^{-((k-1)^2+(k-1)/2)}\left(\zeta_{k,\sigma}^{-1}\right)^{n},
\end{equation*}
where $c_{k,\sigma}$ is some positive constant and the proof of the theorem
is complete. \qed
\section{Combinatorics of ${\sf lv}_k^{\sf 1}$-shapes}\label{S:p'}

\begin{figure}[ht]
\centerline{\includegraphics[width=0.7\textwidth]{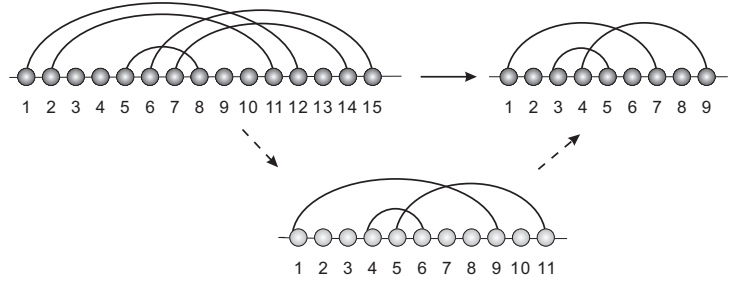}}
\caption{\small
${\sf lv}_k^{\sf 1}$-shapes via the core map and subsequent
identification of unpaired nucleotides: A $3$-noncrossing,
$1$-canonical RNA structure (top-left) is mapped into its ${\sf
lv}_3^{\sf 1}$-shape (top-right).} \label{F:Ik'mapchain}
\end{figure}

\begin{definition}{\bf (${\sf lv}_k^{\sf 1}$-shape)}
Given a $k$-noncrossing, $\sigma$-canonical RNA structure, $\delta$,
its ${\sf lv}_k^{\sf 1}$-shape, ${\sf lv}_k^{\sf 1}(\delta)$, is
derived as follows: first we apply the core map, second we replace a
segment of isolated vertices by a single isolated vertex and third
relabel the vertices of the resulting diagram, see
Fig.\ref{F:Ik'mapchain}.
\end{definition}
More formally, a ${\sf lv}_k^{\sf 1}$-shape is obtained as follows:
if we have a maximal sequence of isolated vertices
$(i,i+1,\ldots,i+\ell')$ (i.e.~ $i-1,i+\ell'+1$ are not isolated),
then we map $(i,i+1,\ldots,i+\ell')\mapsto i$ and if $(i,j)$ is a
arc, it is mapped identically.

Let ${\mathcal C}_k(n,h)$ (${C}_k(n,h)$) denote the set (number) of
$k$-noncrossing core-structures of length $n$ with exactly $h$-arcs.
Let $\mathcal{J}_k(n,h)$ ($j_k(n,h)$) denote the set (number) of
${\sf lv}_k^{\sf 1}$-shapes of length $n$ with $h$-arcs, and let
$j_k(n)$ be the number of all ${\sf lv}_k^{\sf 1}$-shapes of length
$n$ and set
\begin{equation}
{\bf J}_k(z,u)=\sum_{h\geq 0}\sum_{n=2h}^{4h+1}j_k(n,h)z^nu^h \ \text{\rm and}\
{\bf J}_k(z)=\sum_{n\geq 0}j_k(n)z^n.
\end{equation}

\begin{theorem}\label{T:gfIk'}
For $k,n,h\in \mathbb{N}$, $k\geq 2$, the following assertions hold\\
{\rm (a)} the generating functions ${\bf J}_k(z,u)$ and ${\bf J}_k(z)$
are given by
\begin{eqnarray}
\label{E:gfIk'}
{\bf J}_k(z,u) & = & \frac{(1+z)(1+uz^2)}{uz^3+2uz^2+1}
{\bf F}_k\left(\frac{(1+z)^2(1+uz^2)uz^2}
{(uz^3+2uz^2+1)^2}\right)\\
\label{E:o2}
{\bf J}_k(z) & = &\frac{(1+z)(1+z^2)}{z^3+2z^2+1}{\bf F}_k
\left(\frac{(1+z)^2(1+z^2)z^2} {(z^3+2z^2+1)^2}\right).
\end{eqnarray}
{\rm (b)} for $2\leq k\leq 9$, the number of ${\sf lv}_k^{\sf
1}$-shapes of length $n$ satisfies
\begin{equation}\label{E:asy-p'}
j_k(n)\sim c_k'
n^{-((k-1)^2+(k-1)/2)}\left(\mu_{k}'^{-1}\right)^{n},
\end{equation}
where $c_{k}'>0$ and $\mu_{k}'$ is the unique minimum positive real solution of
\begin{equation}\label{E:Iksinner2}
\frac{(1+z)^2(1+z^2)z^2} {(z^3+2z^2+1)^2}=\rho_k^2.
\end{equation}
\end{theorem}
\pff For (a) we consider the map between $k$-noncrossing cores
having exactly $h$ arcs and ${\sf lv}_k^{\sf 1}$-shapes, for $0\leq
h\leq \lfloor\frac{n-1}{2}\rfloor$,
\begin{equation*}
\begin{split}
&\ell\colon {\mathcal C}_k(n,h)\rightarrow\\
&\dot{\bigcup_{b_0\leq b\leq n-2h-1}}\left[ \mathcal{J}_k
(n-b,h)\times \left\{(e_j)_{1 \leq j \leq n-2h-b}\, |\,
\sum_{j=1}^{n-2h-b}e_j=b, e_j \geq 0\right\}\right],
\end{split}
\end{equation*}
where $b_0=\max\{0,n-4h-1\}$. For every
$\beta\in\mathcal{C}_k(n,h)$, $(e_j)_{1\leq j\leq n-2h-b}$ keeps
track of the multiplicities of the deleted isolated vertices. The
map $\ell$ is a (well defined) bijection and
\begin{equation*}
|\{(e_j)_{1 \leq j \leq n-2h-b}\, |\, \sum_{j=1}^{n-2h-b}e_j=b, e_j
\geq 0\}|=\binom{n-2h-1}{b}.
\end{equation*}
We arrive at
\begin{equation*}
C_k(n,h)=\sum_{b=b_0 }^{n-2h-1}\binom{n-2h-1}{b}j_k(n-b,h), \quad
0\leq h\leq \lfloor\frac{n-1}{2}\rfloor.
\end{equation*}
We compute
\begin{equation*}
\begin{split}
\sum_{n\geq 0}\sum_{h=0}^{
\lfloor\frac{n}{2}\rfloor}C_k(n,h)w^{h}x^n=&
\underbrace{
\sum_{n\geq 0}\sum_{h>\lfloor\frac{n-1}{2}\rfloor}^{
\lfloor\frac{n}{2}\rfloor}C_k(n,h)w^{h}x^n}_{{\rm (I)}}+\\
&\underbrace{\sum_{n\geq 0}\sum_{h=0}^{\lfloor\frac{n-1}{2}\rfloor}\sum_{b=b_0
}^{n-2h-1}\binom{n-2h-1}{b}j_k(n-b,h)w^{h}x^n}_{{\rm (II)}},
\end{split}
\end{equation*}
and rewrite (II) as
\begin{equation*}
\begin{split}
&\sum_{n\geq 0}\sum_{h=0}^{\lfloor\frac{n-1}{2}\rfloor}\sum_{b=b_0
}^{n-2h-1}\binom{n-2h-1}{b}j_k(n-b,h)w^hx^n\\
=&\sum_{h\geq0}\sum_{b\geq
0}\sum_{n=2h+b+1}^{4h+b+1}j_k(n-b,h)\binom{n-2h-1}{b}w^hx^n.
\end{split}
\end{equation*}
We derive, setting $s=n-b$,
\begin{equation*}
\begin{split}
=&\sum_{h\geq 0}\sum_{b\geq
0}\sum_{s=2h+1}^{4h+1}j_k(s,h)\binom{s+b-2h-1}{b}w^hx^{s+b}\\
=&\sum_{h\geq 0}\sum_{s=2h+1}^{4h+1}j_k(s,h)\left(\sum_{b\geq
0}\binom{s+b-2h-1}{b}x^b\right)w^hx^{s}\\
=&\sum_{h\geq
0}\sum_{s=2h+1}^{4h+1}j_k(s,h)\left(\frac{x}{1-x}\right)^s
\left((1-x)^2w\right)^h.
\end{split}
\end{equation*}
In view of $j_k(2h,h)=C_k(2h,h)$, we can interpret (I) as follows
\begin{equation*}
\sum_{n\geq
0}\sum_{h>\lfloor\frac{n-1}{2}\rfloor}^{\lfloor\frac{n}{2}\rfloor}
C_k(n,h)w^hx^n=\sum_{h\geq
0}j_k(2h,h)\left(\frac{x}{1-x}\right)^{2h} \left((1-x)^2w\right)^h,
\end{equation*}
which allows for extending the parameter range of $h$
\begin{equation*}
\sum_{n\geq
0}\sum_{h=0}^{\lfloor\frac{n}{2}\rfloor}C_k(n,h)w^hx^n=\sum_{h\geq
0}\sum_{s=2h}^{4h+1}j_k(s,h)\left(\frac{x}{1-x}\right)^s
\left((1-x)^2w\right)^h.
\end{equation*}
Setting $u=(1-x)^2w$ and $z=\frac{x}{1-x}$, we obtain the bivariate
generating function
\begin{equation*}
\sum_{h\geq 0}\sum_{s=2h}^{4h+1}j_k(s,h)z^su^h=\sum_{n\geq
0}\sum_{h=0}^{\lfloor\frac{n}{2}\rfloor}C_k(n,h)\left(u(1+z)^2\right)^h
\left(\frac{z}{1+z}\right)^n.
\end{equation*}
We next consider two power series relations due to
\cite{Reidys:07asym} and \cite{Reidys:07lego}
\begin{eqnarray}
\label{E:Tf}
\sum_{n\geq
0}\sum_{h=0}^{\lfloor\frac{n}{2}\rfloor}\mathrm{T}_{k,1}(n,h)v^{h}y^n
& = &
\frac{1}{vy^2-y+1}
{\bf F}_k\left(\frac{vy^2}{\left(vy^2-y+1\right)^2}\right) \\
\label{E:TC}
\sum_{n\geq 0}\sum_{h=0}^{\lfloor\frac{n}{2}
\rfloor}\mathrm{T}_{k,1}(n,h)v^{h}y^n & = & \sum_{n\geq
0}\sum_{h=0}^{\lfloor\frac{n}{2}\rfloor}C_k(n,h)
\left(\frac{v}{1-vy^2}\right)^hy^n.
\end{eqnarray}
In view of eq.~(\ref{E:Tf}) and eq.~(\ref{E:TC}), we can conclude
\begin{equation*}
\sum_{h\geq 0}\sum_{s=2h}^{4h+1}j_k(s,h)z^su^h=
\frac{(1+z)(1+uz^2)}{uz^3+2uz^2+1}{\bf F}_k\left(\frac{(1+z)^2(1+uz^2)uz^2}
{(uz^3+2uz^2+1)^2}\right)
\end{equation*}
and in particular, setting $u=1$,
\begin{equation*}
{\bf J}_k(z)=\frac{(1+z)(1+z^2)}{z^3+2z^2+1}{\bf F}_k
\left(\frac{(1+z)^2(1+z^2)z^2} {(z^3+2z^2+1)^2}\right),
\end{equation*}
whence assertion (a). \\
Assertion (b) follows in complete analogy to the proof of
Theorem~\ref{T:Ikasy}. First we verify that
the factor $\frac{(1+z)(1+z^2)}{z^3+2z^2+1}$ does not introduce a dominant
singularity of ${\bf J}_k(z)$. Then we verify, using Tab.~\ref{Table:polyroot},
that the unique dominant singularity of ${\bf F}_k\left(\frac{(1+z)^2(1+z^2)z^2}
{(z^3+2z^2+1)^2}\right)$ is the minimum positive real solution of
$\frac{(1+z)^2(1+z^2)z^2} {(z^3+2z^2+1)^2}=\rho_k^2$ for $2\le k\le 9$.
Now (b) follows from
Proposition~\ref{P:algeasym}.
\qed

We finally compute the number of ${\sf lv}^{\sf 1}_k$-shapes induced
by $k$-noncrossing, $\sigma$-canonical RNA structures of fixed
length $n$, ${\sf lv}_{k,\sigma}^{\sf 1}(n)$, setting
\begin{equation}\label{E:lv-3}
{\bf Lv}^{\sf 1}_{k,\sigma}(x)= \sum_{n\geq 0}{\sf
lv}_{k,\sigma}^{\sf 1}(n) x^n.
\end{equation}
\begin{theorem}\label{T:Ik'asy}
Let $k,\sigma\in \mathbb{N}$, where $k\ge 2$. Then the following assertions hold\\
{\rm (a)} the generating function ${\bf Lv}^{\sf 1}_{k,\sigma}(x)$
is given by
\begin{equation}\label{E:oha}
{\bf Lv}^{\sf
1}_{k,\sigma}(x)=\frac{(1+x)(1+x^{2\sigma})}{(1-x)(x^{2\sigma+1}
+2x^{2\sigma}+1)}{\bf F}_k
\left(\frac{(1+x)^2x^{2\sigma}(1+x^{2\sigma})}{\left(x^{2\sigma+1}
+2x^{2\sigma}+1\right)^2}\right)^n.
\end{equation}
{\rm (b)} for $2\leq k\leq 9$ and $1\leq \sigma \leq 10$, we have
\begin{equation}
{\sf lv}^{\sf 1}_{k,\sigma}(n)\sim
c_{k,\sigma}'n^{-((k-1)^2+(k-1)/2)}\left(\chi_{k,\sigma}^{-1}\right)^{n},
\end{equation}
where $c_{k,\sigma}'>0$ and $\chi_{k,\sigma}$ is the unique minimum positive
real solution of
\begin{equation}\label{E:Ik'inner}
\frac{(1+x)^2x^{2\sigma}(1+x^{2\sigma})}{\left(x^{2\sigma+1}
+2x^{2\sigma}+1\right)^2}=\rho_k^2.
\end{equation}
\end{theorem}
\begin{table}[ht]
\begin{center}
\begin{tabular}{cccccccc}
\hline
$\sigma/k$& $2$& $3$ & $4$ & $5$ & $6$ & $7$ & $8$\\
\hline
$1$ &2.09188&4.51263&6.65586&8.73227&10.7804&12.8137&14.8381\\
$2$ &1.56947&2.31767&2.81092&3.21184&3.55939&3.87079&4.15552\\
$3$ &1.38475&1.80408&2.05600&2.24968&2.41081&2.55050&2.67477\\
\hline
\end{tabular}
\centerline{}
\smallskip
\caption{\small The exponential growth rates $\chi_{k,\sigma}^{-1}$
of ${\sf lv}_k^{\sf 1}$-shapes induced by $k$-noncrossing,
$\sigma$-canonical RNA structures of length $n$.} \label{Table:Ik'}
\end{center}
\end{table}
\pff Obviously, we can inflate any structure by adding arcs into its
stacks or duplicating isolated vertices without changing its ${\sf
lv}_k^{\sf 1}$-shape. As a result, we can derive from any ${\sf
lv}_k^{\sf 1}$-shape by inflating its stacks to $\sigma$ arcs, a
unique, minimal, $k$-noncrossing, $\sigma$-canonical structure
inducing it. This observation implies
\begin{equation*}
{\sf lv}_{k,\sigma}^{\sf 1}(n)=\sum_{h=0}^{\lfloor\frac{n}{2\sigma}
\rfloor}\sum_{s=2h}^{\min\{4h+1,n-2(\sigma-1)h\}}j_k(s,h),
\end{equation*}
whence we can rewrite the generating function
\begin{equation*}
{\bf Lv}_{k,\sigma}^{\sf 1}(x)  =  \sum_{h\geq
0}\sum_{s=2h}^{4h+1}\sum_{n\geq 2h(\sigma-1)+s}j_k(s,h)x^{n} =
\frac{1}{1-x}\sum_{h\geq 0}\sum_{s=2h}^{4h+1}
j_k(s,h)x^{2h(\sigma-1)+s}.
\end{equation*}
Employing eq.~(\ref{E:gfIk'}), we derive
\begin{equation*}
{\bf Lv}_{k,\sigma}^{\sf
1}(x)=\frac{(1+x)(1+x^{2\sigma})}{(1-x)(x^{2\sigma+1}
+2x^{2\sigma}+1)}{\bf F}_k
\left(\frac{(1+x)^2x^{2\sigma}(1+x^{2\sigma})}{\left(x^{2\sigma+1}
+2x^{2\sigma}+1\right)^2}\right)^n
\end{equation*}
and assertion (a) follows. As for assertion (b), we proceed in
analogy to the proof of Theorem~\ref{T:Ikasy} and verify that for
$2\leq k\leq 9$ and $1\leq \sigma \leq 10$, the unique minimum
positive real solution, $\chi_{k,\sigma}$, of eq.~(\ref{E:Ik'inner})
is the unique dominant singularity of generating function ${\bf
Lv}_{k,\sigma}^{\sf 1}(x)$. Consequently, Proposition
\ref{P:algeasym} implies that
\begin{equation*}
{\sf lv}_{k,\sigma}^{\sf 1}(n)\sim
c_{k,\sigma}'n^{-((k-1)^2+(k-1)/2)}\left(\chi_{k,\sigma}^{-1}\right)^{n},
\end{equation*}
where $c_{k,\sigma}'$ is some positive constant, whence (b) and the theorem
is proved. \qed

\section{Conclusion}
${\sf lv}_{k}^{\sf 1}$- and ${\sf lv}_{k}^{\sf 5}$-shapes of $k$-noncrossing,
$\sigma$-canonical RNA pseudoknot structures provide a significant simplification
of complicated molecular configurations with cross-serial interactions.
The asymptotic formulas presented in Theorem~\ref{T:Ikasy} and
Theorem~\ref{T:Ik'asy}
\begin{eqnarray*}
{\sf lv}_{k,\sigma}^{\sf 5}(n) &\sim &
c_{k,\sigma}n^{-((k-1)^2+(k-1)/2)}\left(\zeta_{k,\sigma}^{-1}\right)^{n}\\
{\sf lv}^{\sf 1}_{k,\sigma}(n) & \sim &
c_{k,\sigma}'n^{-((k-1)^2+(k-1)/2)}\left(\chi_{k,\sigma}^{-1}\right)^{n},
\end{eqnarray*}
imply all asymptotic results on abstract shapes of secondary structures
in the literature (note $n^{-((k-1)^2+(k-1)/2)}=n^{-3/2}$).

The growth rates of ${\sf lv}^{\sf 1}_{k}$- and ${\sf lv}_{k}^{\sf 5}$-shapes
of $k$-noncrossing, $\sigma$-canonical structures, are displayed in
Tab.~\ref{Table:Tk2} and Tab.~\ref{Table:Tk3}, where they are contrasted with
the exponential growth rates of $k$-noncrossing, $\sigma$-canonical
structures, $\gamma_{k,\sigma}$.
\begin{table}[ht]
\begin{center}
\begin{tabular}{cccccccc}
\hline
$k$& $2$ & $3$ & $4$ & $5$ & $6$ & $7$ & $8$ \\
\hline
${\gamma}_{k,2}^{-1}$  &1.96798&2.58808&3.03825&3.41383&3.74381&4.04195&4.31617\\
$\chi_{k,2}^{-1}$ &1.56947&2.31767&2.81092&3.21184&3.55939&3.87079&4.15552\\
$\zeta_{k,2}^{-1}$ &1.26585&1.93496&2.41152&2.80275&3.14338&3.44943&3.72983\\
\hline
\end{tabular}
\centerline{}
\smallskip
\caption{\small The exponential growth rates of arbitrary $k$-noncrossing, $2$-canonical
RNA structures of length $n$ and the numbers of their induced
${\sf lv}_k^{\sf 1}$ and
${\sf lv}_k^{\sf 5}$ shapes.} \label{Table:Tk2}
\end{center}
\end{table}

\begin{table}[ht]
\begin{center}
\begin{tabular}{cccccccc}
\hline
$k$& $2$ & $3$ & $4$ & $5$ & $6$ & $7$ & $8$ \\
\hline
$\gamma_{k,3}^{-1}$  & 1.71599&2.04771&2.27036&2.44664&2.59554&2.72590&2.84267 \\
$\chi_{k,3}^{-1}$ &1.38475&1.80408&2.05600&2.24968&2.41081&2.55050&2.67477\\
$\zeta_{k,3}^{-1}$ & 1.17928& 1.55752& 1.80082 &1.98945&2.14693&2.28376&2.40567 \\
\hline
\end{tabular}
\centerline{}
\smallskip
\caption{\small The exponential growth rates of arbitrary $k$-noncrossing, $3$-canonical
RNA structures of length $n$ and the numbers of their induced
${\sf lv}_k^{\sf 1}$ and ${\sf lv}_k^{\sf 5}$ shapes.} \label{Table:Tk3}
\end{center}
\end{table}
Table~\ref{Table:Tk3} shows that the exponential growth rate of
${\sf lv}_{3}^{\sf 5}$-shapes of $k$-noncrossing $3$-canonical
structures are significantly smaller than that of all
$k$-noncrossing $3$-canonical structures. Therefore, the abstract
${\sf lv}_{3}^{\sf 5}$-shapes represent a meaningful reduction.
At {\tt http://www.combinatorics.cn/cbpc/paper.html}, we provide supplemental material for our results.

{\bf Acknowledgments.}
This work was supported by the 973 Project, the PCSIRT Project of
the Ministry of Education, the Ministry of Science and Technology,
and the National Science Foundation of China.

\bibliographystyle{plain}


\end{document}